\documentclass[12pt]{article}
\usepackage{amsmath, amsthm, amsfonts}
\usepackage{amssymb}
\usepackage{amscd}
\usepackage{verbatim}

\begin{document}
\newcommand{\core}{{C_c^\infty(\R^N)}}
\newcommand{\crit}{{\frac{2N}{N-2}}}
\newcommand{\dom}{{\R^N}}
\newcommand{\bdry}{{\R^N}}
\newcommand{\enn}{N}
\newcommand{\merhav}{{\mathcal D}^{1,2}}
\newcommand{\be}{\begin{equation}}
\newcommand{\ee}{\end{equation}}
\newcommand{\bea}{\begin{eqnarray}}
\newcommand{\eea}{\end{eqnarray}}
\newcommand{\bean}{\begin{eqnarray*}}
\newcommand{\eean}{\end{eqnarray*}}
\newcommand{\thkl}{\rule[-.5mm]{.3mm}{3mm}}
\newcommand{\cw}{\stackrel{D}{\rightharpoonup}}
\newcommand{\id}{\operatorname{id}}
\newcommand{\supp}{\operatorname{supp}}
\newcommand{\wlim}{\mbox{ w-lim }}
\newcommand{\mymu}{{x_N^{-p_*}}}
\newcommand{\R}{{\mathbb R}}
\newcommand{\N}{{\mathbb N}}
\newcommand{\Z}{{\mathbb Z}}
\newcommand{\Q}{{\mathbb Q}}
\newtheorem{theorem}{Theorem}[section]
\newtheorem{corollary}[theorem]{Corollary}
\newtheorem{lemma}[theorem]{Lemma}
\newtheorem{definition}[theorem]{Definition}
\newtheorem{remark}[theorem]{Remark}
\newtheorem{proposition}[theorem]{Proposition}
\newtheorem{conjecture}[theorem]{Conjecture}
\newtheorem{question}[theorem]{Question}
\newtheorem{example}[theorem]{Example}
\newtheorem{Thm}[theorem]{Theorem}
\newtheorem{Lem}[theorem]{Lemma}
\newtheorem{Pro}[theorem]{Proposition}
\newtheorem{Def}[theorem]{Definition}
\newtheorem{Exa}[theorem]{Example}
\newtheorem{Exs}[theorem]{Examples}
\newtheorem{Rems}[theorem]{Remarks}
\newtheorem{Rem}[theorem]{Remark}
\newtheorem{Cor}[theorem]{Corollary}
\newtheorem{Conj}[theorem]{Conjecture}
\newtheorem{Prob}[theorem]{Problem}
\newtheorem{Ques}[theorem]{Question}
\newcommand{\pf}{\noindent \mbox{{\bf Proof}: }}


\renewcommand{\theequation}{\thesection.\arabic{equation}}
\catcode`@=11 \@addtoreset{equation}{section} \catcode`@=12


\title{Concentration-compactness principle for mountain pass problems}
\author{Kyril Tintarev\thanks{Research done while visiting Technion -- Israel Institute of Technology.
Supported in part by a grant from Swedish Research Council.}
\\{\small Department of Mathematics}\\{\small Uppsala University}\\
{\small SE-751 06 Uppsala, Sweden}\\{\small
kyril.tintarev@math.uu.se}}
\maketitle
\newcommand{\dnorm}[1]{\thkl #1 \thkl\,}

\begin{abstract}
In the paper we show that critical sequences associated with the
mountain pass level for semilinear elliptic problems on $\R^N$
converge when the non-linearity is subcritical, superlinear and
satisfies the penalty condition $F_\infty(s)<F(x,s)$. The proof
suggests a concentration compactness framework for minimax
problems, similar to that of P.-L.Lions for constrained
minima.\vskip3mm
\noindent  2000 {\em Mathematics Subject Classification:} 35J20, 35J60, 49J35
\end{abstract}


\section{Introduction}
In this paper we prove an existence result for the classical
semilinear elliptic problem on $\R^N$: \be \label{Theeq} -\Delta
u+u= f(x,u), u\in H^1(\R^N).\ee

The classical existence proof for the analogous Dirichlet problem
on bounded domain, based on the mountain pass lemma of \cite{AR}
fails in the case of $\R^N$, since the Palais-Smale condition does
not anymore follow from compactness of Sobolev imbeddings and a
concentration compactness argument is needed. There are numerous
publications where the concentration compactness is used in
minimax problems, including the problem considered here (a
representative bibliography on the subject can be found in the
books of Chabrowski \cite{Chab} and Willem \cite{Willem}), but
given that in problems on $\R^N$ the $(PS)_c$ condition fails,
typically, for every $c$ that is is a linear combination, with
positive integer coefficients, of critical values, the
Palais-Smale condition has been proved only with severe
restrictions on the nonlinearity $f(x,s)$.

We consider here a set of conditions on the functional, similar to
the conconcentration compactness framework as set by P.-L.Lions
(\cite{PLL1a},\cite{PLL1b},\cite{PLL2a}, \cite{PLL2b}), where
conditions for existence of minima can be formulated as the
following prototype assumptions: a) the functionals are
continuous, b) critical sequences are bounded (achieved by
regarding constrained minima), c) constrained minimal values are
subadditive with respect to the parameter of constraint level, and
d) the functionals are invariant relative to transformations
causing the loss of compactness or their asymptotic values (with
respect to the unbounded sequences of the transformations) satisfy
a penalty condition.

In the present paper we consider a functional $G$ on a Hilbert
space $H$ with an asymptotic (with respect to unbounded sequence
of transformations responsible for loss of compactness) value
$G_\infty$. Let $\Phi$ be an appropriate set of mappings from a
metric space $X$ into $H$ with fixed values on a $X_0\subset X$,
and let $\rho:=\sup G(\varphi(X_0)<c=:\inf_{\varphi\in\Phi}\sup
G(\varphi(X))$. We regard the following heuristic conditions,
whose formal counterparts for the functional associated with
(\ref{Theeq}) will are given in Section 2.

 $a'$) $G\in C^1(H)$ (in the semilinear elliptic
 case, subcritical growth of $f$);

 $b'$) critical sequences at the level $c$ are bounded (in the semilinear elliptic
 case follows from an assumption of superlinearity for $f$);

 $c'$) all critical points of $G_\infty$ with critical values in $(\rho,c]$
have the Morse index greater or equal to the one associated with
the minimax (a weaker version: for every critical point $w$ of
$G_\infty$ such that $\rho<G_\infty(w)\le c$ there is a sequence
of paths $\varphi_k$ such that $d(w,\varphi_k(X))\to 0$ and
$\sup_{x\in X}G_\infty(\varphi_k(x))\to G(w)$) - in the semilinear
elliptic case with a mountain pass, the sufficient condition is
$s\mapsto f_\infty(s)/|s|$ monotone increasing; and

$d'$) invariance or penalty condition ($f(x,s)=f(s)$ or
$F(x,s)>F_\infty(s)$ where $F(x,t)=\int_0^tf(x,s)ds$).

Existence of critical points is proved by verifying $(PS)_c$ for a
single value $c$, namely the one given by the mountain pass
statement. Sharp estimates of $c$ are based on the global
compactness theorem by I.Schndler and the author (\cite{acc}),
which is a functional-analytic generalization of earlier
"multi-bump" weak convergence lemmas (Struwe,\cite{Struwe}; Lions,
\cite{Lions86}; Cao and Peng \cite{Cao}.

\section{Existence theorem}
We consider the Hilbert space $H^1(\R^N)$, $N\in\N$, defined as
the completion of $\core$ with respect to the norm \be
\|u\|^2=\int_{\R^N}|\nabla u|^2+u^2. \ee In what follows the
notation of norm without other specification will refer to this
$H^1$-norm. The space $H^1(\R^N)$ is continuously imbedded into
$L^p(\R^N)$ for $2\le p\le \frac{2N}{N-2}$ when $N>2$ and for
$p\ge 2$ for $N=1,2$. For convenience we set $2^*=\frac{2N}{N-2}$
for $N>2$ and $2^*=\infty$ for $N=1,2$. Let $f:\R^N\times\R$ be
continuous function and let \be F(x,t)=\int_0^tf(x,s)ds,\ee \be
g(u)=\int_{\R^N} F(x,u(x))dx,\ee

and

\be \label{G} G(u)=\frac12\|u\|^2-g(u). \ee

We assume that $f(x,s)\to f_\infty(s)$ as $|x|\to\infty$ and
follow the definitions above to define $F_\infty$, $g_\infty$ and
$G_\infty$.

Let $\omega\subset H^1(\R^N)$ be a path-connected component of
$G^{-1}(-\infty,-1)$ containing infinity, that is, the union of
all path-connected subsets $\omega^\prime\subset
G^{-1}(-\infty,-1)$, such that, for any $R>0$ every point of
$\omega^\prime$ can be connected by an arc in $G^{-1}(-\infty,-1)$
to a point $e$ with $\|e\|>R$. Let

\be \Phi=\{\varphi\in C([0,1];H^1(\R^N)): \varphi(0)=0,
\varphi(1)\in \omega \}. \ee

\begin{theorem}
\label{main} Assume that:

for every $\epsilon>0$ there exist $p_\epsilon\in(2,2^*)$ and
$C_\epsilon>0$ such that
$$
\leqno{(A)}\; |f(x,s)|\le
\epsilon(|s|+|s|^{2^*-1})+C_\epsilon|s|^{p_\epsilon-1}, s\in\R,
x\in\R^N;
$$

There exists a $\mu>2$, such that
$$
\leqno{(B)}\; f(x,s)s\ge\mu F(x,s), s\in\R, x\in\R^N;
$$

$$
\leqno{(C)}\; s\mapsto f_\infty(s)/|s|, s\in\R, \mbox{ is
increasing};
$$

$$
\leqno{(D)}\; F(x,s)>F_\infty(s), s\in\R\setminus\{0\}, x\in\R^N.
$$
Then $\Phi\neq\emptyset$; \be \label{c}
c:=\inf_{\varphi\in\Phi}\max_{t\in[0,\infty)} G(\varphi(t))>0;\ee
there is a sequence $u_k\in H^1(\R^N)$ such that $G^\prime(u_k)\to
0$, $G(u_k)\to c$; every such sequence has a subsequence
convergent in $H^1(\R^N)$. Consequently, $u=\lim u_k$ satisfies
$G(u)=c$ and $G^\prime(u)=0$ (and therefore, $u$ is a solution of
(\ref{Theeq})).
\end{theorem}

Condition (A) is a well-known sufficient condition for $G\in
C^1(H^1(\R^N ))$.

\begin{lemma}
\label{lem} Let $G$ be as in (\ref{G}). Assume conditions (A) and
(C) of Theorem~\ref{main}. Then for every $w\in
H^1(\R^N)\setminus\{0\}$, the path $\varphi(t)=tw$,
$t\in(0,\infty)$, is in $\Phi$ and the constant (\ref{c}) is
positive. If, in addition, $G^\prime(w)=0$, then $\max_t
G_\infty(\varphi(t))$ is attained at $\varphi(1)=w$.
\end{lemma}

\begin{proof} The first assertion of the lemma follows easily
from (C) and the second is a consequence of (A) (the proof is a
trivial modification of the one in \cite{AR}).

Let $w\neq 0$ satisfy $G^\prime(w)=0$. From (C) follows that the
function $s\mapsto s^{-1}\frac{d}{ds}G_\infty(sw^{(n)})$ is
decreasing on $(0,\infty)$. Then, since

\begin{eqnarray*}
\frac{d}{ds}G_\infty(sw)=s\|w\|^2-\int f_\infty(sw)w\\
=s\left(\|w\|^2-\int \frac{f_\infty(sw)}{sw}w^2dx\right),
\end{eqnarray*}
the function $s\mapsto\gamma(s):= G_\infty(sw)$ has at most one
critical point. Since $\gamma(0)=0$, $\gamma(s)<0$ for $s$ large
and has positive values (because $c>0$), the critical point of
$\gamma$ is a point of maximum. Since $G^\prime(w)=0$,
$(G^\prime(w),w)=0$, which is equivalent to $\gamma^\prime(1)=0$.
Since $\gamma(s)$ has a unique critical point, which is a point of
maximum, $s\mapsto G(sw)$ attains its maximum at $s=1$.
\end{proof}


\section{Global compactness}
In this section we present statements from \cite{acc} concerning
weak convergence that will be used in the proof.

\begin{theorem}
\label{acc}
Let $u_{k}  \in H$ be a  bounded sequence. Then there
exists $w^{(n)}  \in H$, $g_{k} ^{(n)}  \in D$, $k,n \in
\mathbb{N}$ such that  for a renumbered subsequence
\begin{gather}
\label{separates} g_k^{(1)}=id,\; {g_{k} ^{(n)}} ^{-1}  g_{k}
^{(m)}\rightharpoonup 0 \mbox{ for } n \neq m,
\\
w^{(n)}=\wlim  {g_{k} ^{(n)}}^{-1}u_k
\\
\label{norms} \sum_{n \in \mathbb{N}} \|w ^{(n)}\|^2 \le \limsup
\|u_k\|^2
\\
\label{BBasymptotics} u_{k}  - \sum_{n\in\mathbb{N}}  g_{k} ^{(n)}
w^{(n)}  \cw 0.
\end{gather}
\end{theorem}

In particular $u\mapsto u(\cdot-y)$, $y\in\Z^N$, form a
dislocation group in $H^1_0(\R^N)$, and $u_k\cw 0$ is equivalent
to $u_k\to 0$ in $L^p(\R^N)$, $p\in(2,2^*)$ (an equivalent
statement is found in \cite{Lieb}).

The following lemma is similar to the Brezis-Lieb lemma from
\cite{BrLieb} and is a trivial modification of analogous lemma
from \cite{acc}.

\begin{lemma}
\label{BL} Assume that $F$ satisfies (A) and that $u_k$ and
$(w^{(n)})$ are as in Theorem~\ref{acc}. Then \be \label{EqBL}
\int F(x,u_k)\to \int F(x,w^{(1)})+\sum_{n\ge 2}\int F_\infty
(w^{(n)}). \ee
\end{lemma}

\section{Proof of Theorem~\ref{main}}
Step 1. By Lemma~\ref{lem}, $\Phi\neq\emptyset$ and $c>0$. By (A),
(\ref{c}) and the mountain pass lemma (\cite{AR}), there is a
sequence $u_k$ such that $G^\prime(u_k)\to 0$ and $G(u_k)\to c$.
By (B), $u_k$ is bounded in $H^1$ (see, again, the argument of
\cite{AR}) and we can apply Theorem~\ref{acc}, referring in what
follows to the renamed subsequence. By (\ref{EqBL}) and
(\ref{norms}),

\begin{equation}
\label{frombelow1} c\ge \frac12\sum_{n\in\N}\|w^{(n)}\|^2-\int
F(x,w^{(1)})- \sum_{n\ge 2}\int F_\infty(w^{(n)}).
\end{equation}

From $G^\prime(u_k)\to 0$ (since (A), by compactness of local
imbeddings of $H^1$ into $L^p$ implies weak convergence of
$g^\prime(u_k)$ for weakly convergent $u_k$) follows

\begin{equation}
\label{crbalance} \|w^{(1)}\|^2=\int
f(x,w^{(1)})w^{(1)}\;,\|w^{(n)}\|^2=\int f_\infty(w^{(n)})w^{(n)}
\mbox{ for } n\ge 2.
\end{equation}

Substituting (\ref{crbalance}) into (\ref{frombelow1}) we get

\begin{equation}
\label{frombelow} c\ge \int
\left(\frac12f(x,w^{(1)})w^{(1)}-F(x,w^{(1)})\right)+\sum_{n\ge
2}\int \left(\frac12
f_\infty(w^{(n)})w^{(n)}-F_\infty(w^{(n)})\right).
\end{equation}

Step 2. Note that
\begin{equation}
\label{positive tilda} \frac12f(x,s)s-F(x,s)>0\mbox{ and
}\frac12f_\infty(s)s-F_\infty(s)>0, s\neq 0.
\end{equation}
The first relation follows from (B): \be \label{tt}
\frac12f(x,s)s-F(x,s)\ge (\frac{\mu}{2}-1)F(x,s),\ee
$F(x,s)>F_\infty(s)$ by (D) and $F_\infty(s)>0$ for $s\neq 0$ due
to (A). The second relation follows going to the limit in
(\ref{tt}) as $|x|\to\infty$ and using positivity of
$F_\infty(s)$.

Step 3. Assume that

\be \label{false} w^{(n)}\neq 0\mbox{ for some } n\neq 1.\ee Let
us estimate $c$ from above by choosing paths $s\mapsto
sw^{(n)}(\cdot-y_k)\in H^1(\R^N)$ with $y_k\in\Z^N$,
$|y_k|\to\infty$. Then
\begin{equation}
c\le \sup_{s\in(0,\infty)} G(sw^{(n)}(\cdot-y_k)).
\end{equation}
By taking $k\to\infty$ we have
\begin{equation}
\label{straightpath} c\le \sup_{s\in(0,\infty)}
G_\infty(sw^{(n)}).
\end{equation}
By Lemma~\ref{lem},
$\sup_{s\in(0,\infty)}G_\infty(sw^{(n)})=G_\infty(w^{(n)})$, and,
therefore,

\be c\le G_\infty(w^{(n)}). \ee

Comparing this with (\ref{frombelow}), we see, due to
(\ref{positive tilda}), that for $m\neq n$, $w^{(m)}=0$ with
necessity and therefore

\begin{equation}
\label{wrong-c} c=G_\infty(w^{(n)}).
\end{equation}

This is clearly false: consider a path $s\mapsto sw^{(n)}$. Then
by (D), $\sup_{s}G(sw^{(n)})<\sup_{s}G_\infty
(sw^{(n)})=G_\infty(w^{(n)})=c$, which contradicts the definition
of $c$.

We conclude that the assumption (\ref{false}) is false and
$w^{(n)}=0$ for all $n\neq 1$.

Step 4. We conclude from Step 4 and (\ref{BBasymptotics}) that
$u_k\to w^{(1)}$ in $L^r$ for any $r\in (2,2^*)$. Then from (A)
follows $g^\prime(u_k)\to g^\prime(w^{(1)})$, and, since
$u_k-g^\prime(u_k)\to 0$, $u_k$ is a convergent sequence in
$H^1(\R^N)$. We conclude that $u_k\to w^{(1)}$ in $H^1(\R^N)$. By
continuity, $G^\prime(w^{(1)})=0$ and $G(w^{(1)})=c$.

\end{document}